\documentclass[Conference]{IEEEtran}
\usepackage{multirow}
\usepackage{multicol}
\usepackage{lipsum}
\usepackage{graphicx}
\usepackage{graphics}
\usepackage{amsmath}
\usepackage{amsbsy}
\usepackage{blindtext}
\usepackage{tcolorbox}
\usepackage{amssymb}
\usepackage{scalerel}
\usepackage{mathabx}
\usepackage{verbatim}
\usepackage{booktabs}
\usepackage{rotating,tabularx}
\usepackage{mathrsfs} 
\usepackage{hyperref}

\usepackage{adjustbox}
\usepackage{soul}
\usepackage{xcolor}
\usepackage{cite}
\usepackage{tikz}
\usepackage{color, colortbl}
\usepackage[first=0,last=9]{lcg}
\definecolor{LightGray}{rgb}{0.7,0.7,0.7}

\usepackage[wby]{callouts}
\usetikzlibrary{shapes.multipart}

\usepackage{array}
\usepackage{makecell}

\allowdisplaybreaks

\usepackage{amsthm}
\theoremstyle{definition}

\theoremstyle{remark}

\usepackage[utf8]{inputenc}
\usepackage[english]{babel}

\usepackage[caption=false,font=footnotesize]{subfig}
\usepackage[T1]{fontenc}
\usepackage{scalerel,stackengine}
\newcommand\reallywidecheck[1]{%
\savestack{\tmpbox}{\stretchto{%
  \scaleto{%
    \scalerel*[\widthof{\ensuremath{#1}}]{\kern-.6pt\bigwedge\kern-.6pt}%
    {\rule[-\textheight/2]{1ex}{\textheight}}
  }{\textheight}%
}{0.5ex}}%
\stackon[1pt]{#1}{\scalebox{-1}{\tmpbox}}%
}
\stackMath
\IEEEoverridecommandlockouts
\IEEEoverridecommandlockouts

\newcommand*{\mrn}{\textcolor{black}}

\newif\ifarxiv
\arxivtrue
\arxivfalse

\newcolumntype{s}{>{} p{3.0cm}}
\newcolumntype{b}{>{} p{4.0cm}}

\renewcommand{\baselinestretch}{.94}

\usepackage[margin=1in,footskip=0.4in]{geometry}

\begin{document}

\title{\LARGE\bf
Investigating the Impact of Electric Vehicle Charging Loads on CSUN's Electric Grid}

\author{Daniel Garcia Aguilar, Logan DeHay, Jahn Aquino, Erik Jensen, Juan Rodriguez, Mohammad Rasoul Narimani, Silvia Carpitella, Kourosh Sedghisigarchi, and Xudong Jia
}

\maketitle

\begin{abstract}
This paper examines the impact of electric vehicle (EV) charging stations on the capacity of distribution feeders and transformers within the electric grid at California State University Northridge (CSUN). With the increasing adoption of both residential and commercial EVs and the rapid expansion of EV charging infrastructure, it is critical to evaluate the potential overloading effects of intensive EV charging on power distribution systems. This research assesses the impact of EV charging on the operation of CSUN's electric grid, identifying potential overload risks under projected EV adoption scenarios. Detailed simulations and analyses are conducted to quantify the extent of these impacts, focusing on various levels of EV penetration and charging patterns. The study also explores the impact of distributed generation on reducing the stress incurred by EV loads. The findings provide essential insights for utility companies, highlighting the need for strategic upgrades to distribution systems. These insights will help in developing robust strategies for both current operations and future planning to accommodate growing EV charging demands, ensuring grid stability and reliability in the face of increasing electrification of the transportation sector. The modeled CSUN electric grid can be used as a benchmark to study the impact of EV loads in dense areas on various parameters of the grid.
\end{abstract}

\section{Introduction}
\label{Introduction}

The rapid growth in electric vehicle (EV) usage is significantly increasing the demand for electricity needed to charge these vehicles. Lower prices for batteries and government incentives fuel this trend, leading to the production of more EVs, which in turn increases electricity demand~\cite{muratori2021rise}. This surge in demand places a considerable strain on the electric power grid, which has limited capacity for transmission and distribution. As a result, the rise in electricity demand from EVs can lead to issues such as line congestion, highlighting the need for grid upgrades. Therefore, it is crucial to investigate the impact of EV load on the electric grid to make informed decisions, ensuring the grid can handle the additional load without compromising service quality~\cite{secchi2023smart}. 

Extensive research has been dedicated to mitigating the impact of EV load on the power grid. One promising solution is smart charging, which intelligently manages the timing of vehicle charging by considering factors such as renewable energy availability and grid load. The approach in~\cite{crozier2020opportunity} examines where and when vehicles are most commonly used, their power consumption, and the structure of the electricity network. Using conditional probability modeling, it simulates various charging scenarios to ensure a smooth transition to electric vehicles without overwhelming the power grid infrastructure~\cite{crozier2020opportunity}.

In managing the charging of EVs, two primary strategies are employed: centralized and decentralized~\cite{yong2023electric, narimani2019efficient}. In the centralized approach, a single entity, such as a utility company or charging network operator, schedules when groups of EVs will charge. Conversely, the decentralized method allows individual EVs to decide when to charge based on signals they receive regarding electricity prices. Both strategies aim to achieve several key objectives: smoothing out electricity demand across different times and locations~\cite{ramadan2018smart}, minimizing energy loss during charging~\cite{nafisi2015two}, ensuring the balance of the electrical system's phases, and promoting the use of renewable energy sources~\cite{clairand2017tariff, asghari2019method}. Most charging algorithms are designed to ensure that EVs reach their target charge levels efficiently without exceeding predetermined energy limits~\cite{yuvaraj2024comprehensive}. Some algorithms add constraints like monitoring the electrical network's performance to optimize the charging process and improve system efficiency~\cite{yuvaraj2024comprehensive, nour2019smart,narimani2019demand, narimani2017multi,narimani2015dynamic}.

The power distribution systems are designed with a restricted overload threshold, meaning that increased EV charging loads would necessitate either fortifying the current infrastructure or introducing new capacities~\cite{yuvaraj2024comprehensive,azizivahed2017new}. Different charging stations, slow and fast charging, have various charging levels and incur different load levels on the power grid, and thus their impact on the power grid is different. One solution to alleviate the impact of charging stations, specifically fast charging stations, is by linking them to medium transmission lines. This would minimally impact line characteristics due to the higher voltage levels and relatively smaller power output of EV charging stations compared to the existing transmission system load~\cite{yuvaraj2024comprehensive}. A comprehensive review of the impact of EV charging stations in the power grid can be found in~\cite{rahman2022comprehensive}.

In dense communities like university campuses, where numerous students and staff want to charge their vehicles upon arrival, the demand for electricity peaks significantly. If no immediate action is taken to ameliorate this condition, it might damage the system or trigger overload relays, leading to load losses and failure propagation~\cite{boyaci2022spatio, narimani2016reliability}. Thus, it is important to ensure that the electric grid can withstand the electricity peak caused by plugging numerous EVs into the grid.  In this connection, this paper investigates the impact of actual EV loads on California State University Northridge's (CSUN's) electric power grid. 

To analyze the EV's load impact on the grid, we first model CSUN's electric power grid using PowerWorld software.
We gathered load data from the installed meters in each building and charging station on the campus every 15 minutes. Next, we identified the transformer and cable types, as well as the length of each cable, to calculate the impedance of these components. To analyze the impact of load variation on the grid, we leveraged the PowerWorld interface, Automation Server (SimAutho)~\cite{thayer2020easy}, which connects the PowerWorld environment with any programming language. By writing a Python script, we automated the process of changing various quantities, including the load at buses, in our model and executed the power flow analysis for various load profiles. Our models represent a realistic simulation of the electric power distribution grid, capable of estimating the impact of EVs on power distribution systems by analyzing line overloads and voltage levels. With electric distribution system models extending down to individual buildings, the proposed approach offers a more granular level of detail than what is typically possible with analytical methods. Thus, we envision this model serving as a validation platform for tasks such as identifying relevant modeling approximations and evaluating control strategies developed with less detailed models. The primary contribution of this paper is the presentation of a ready-to-use model for detailed simulations of EV charging operations in small-size power distribution systems.

\begin{figure*}[t]
 \centering
\captionsetup{justification=centering}
\begin{tikzpicture}
\node at (0,0) {\includegraphics[scale=1.45,trim= 7cm 7.2cm 10cm 7.15cm  2.3cm,clip]{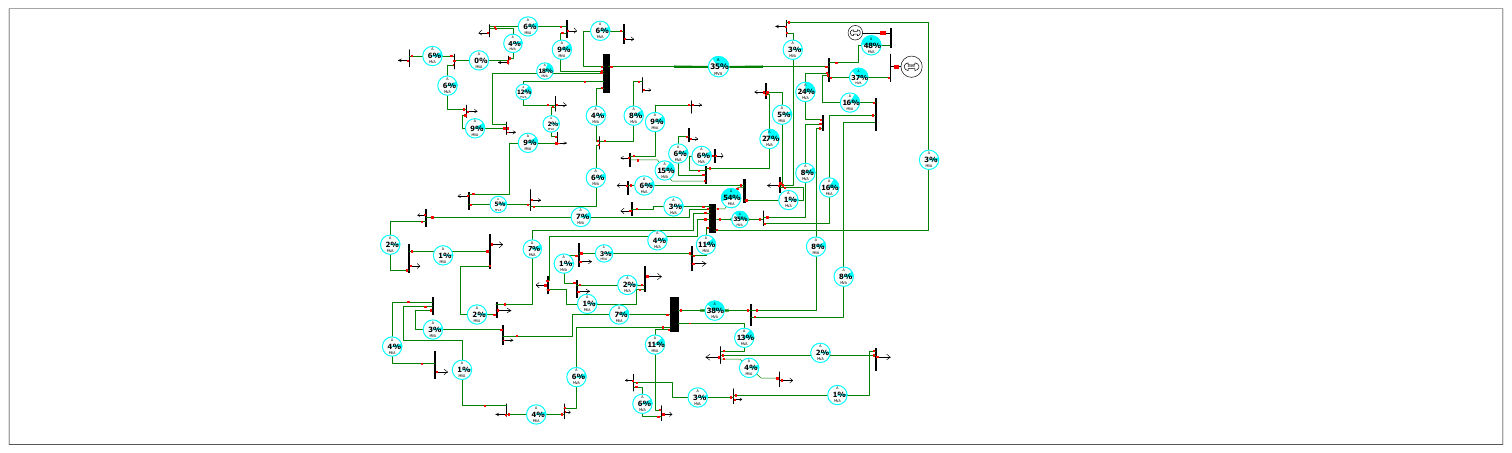}};
 
\node [anchor=west, font=\scriptsize] at (-5.5,5.1) {Physical Plant};
\node [anchor=west, font=\scriptsize] at (-7.5,4.5) {Extended Learning};
\node [anchor=west, font=\scriptsize] at (-4.6,3.8) {Juniper Hall};
\node [anchor=west, font=\scriptsize] at (-3.6,5.2) {Art Design Center};
\node [anchor=west, font=\scriptsize] at (-1,4.8) {E6 Mathador Hall};
\node [anchor=west, font=\scriptsize] at (-1.5,4) { SubB (LV)};
\node [anchor=west, font=\scriptsize] at (-6.2,2.7) {Education};
\node [anchor=west, font=\scriptsize] at (-4.6,2.6) {Jacaranda Hall};
\node [anchor=west, font=\scriptsize] at (-3.4,3.3) {Redwood Hall};
\node [anchor=west, font=\scriptsize] at (-0.7,3.5) {Sustainability Center};
\node [anchor=west, font=\scriptsize] at (0.5,2.8) {SU(AN)};
\node [anchor=west, font=\scriptsize] at (1.8,3.6) {SU(CP)};
\node [anchor=west, font=\scriptsize] at (2.2,5.2) {Athletic Field};
\node [anchor=west, font=\scriptsize] at (3.2,4.2) {SubB};
\node [anchor=west, font=\scriptsize] at (4.8,3.4) {DWP Pole};
\node [anchor=west, font=\scriptsize] at (4.5,5.1) {Tampa DWP Pole};
\node [anchor=west, font=\scriptsize] at (0.5,2.3) {SU Center};
\node [anchor=west, font=\scriptsize] at (1.1,1.8) {SU (AE)};
\node [anchor=west, font=\scriptsize] at (1.1,1.25) {Student REC};
\node [anchor=west, font=\scriptsize] at (1.9,0.4) {Sub A};
\node [anchor=west, font=\scriptsize] at (0.05,0.75) {Sub A (LV)};
\node [anchor=west, font=\scriptsize] at (-1.6,0.6) {Street LTS};
\node [anchor=west, font=\scriptsize] at (-1.6,1.1) {SU (As)};
\node [anchor=west, font=\scriptsize] at (-1.6,1.9) {SU Admin};
\node [anchor=west, font=\scriptsize] at (-3.6,1.7) {SEQUOIA Hall};
\node [anchor=west, font=\scriptsize] at (-3.9,1) {Oviatt library};
\node [anchor=west, font=\scriptsize] at (-5.5,1) {Byramian};
\node [anchor=west, font=\scriptsize] at (-6.9,0.5) {University Hall};
\node [anchor=west, font=\scriptsize] at (-6.7,-0.3) {Seirra Center};
\node [anchor=west, font=\scriptsize] at (-5.2,-0.1) {Jerome Richfield};
\node [anchor=west, font=\scriptsize] at (-6.8,-1.7) {Parking B2};
\node [anchor=west, font=\scriptsize] at (-4.8,-1.7) {Seirra Hall};
\node [anchor=west, font=\scriptsize] at (-4.8,-3) {Manzanita Hall};
\node [anchor=west, font=\scriptsize] at (-6.8,-3.9) {Parking B3};
\node [anchor=west, font=\scriptsize] at (-4.9,-4.9) {Nordhoff Hall};
\node [anchor=west, font=\scriptsize] at (-3.2,-4.2) {Cypress Hall};
\node [anchor=west, font=\scriptsize] at (-1.5,-3.5) {Bookstore};
\node [anchor=west, font=\scriptsize] at (-0.5,-4.9) {Soraya Hall};
\node [anchor=west, font=\scriptsize] at (1,-4.5) {Monterey Hall};
\node [anchor=west, font=\scriptsize] at (0,-2.9) {Health Center};
\node [anchor=west, font=\scriptsize] at (1.3,-1.8) {Sub C (HV)};
\node [anchor=west, font=\scriptsize] at (-0.3,-1.6) {Sub C};
\node [anchor=west, font=\scriptsize] at (-0.1,-1.2) {Sattelite Plant};
\node [anchor=west, font=\scriptsize] at (-1.5,-0.9){Chapparal Hall};
\node [anchor=west, font=\scriptsize] at (2.3,-3.5){Chrisholm Hall};
\node [anchor=west, font=\scriptsize] at (5.1,-3.5){Parking G3};
\end{tikzpicture}%
 \caption{CSUN's campus electric grid with regular load, excluding projected EV loads, shows no line congestion.}
	\label{fig:powerworld_CSUN}
\end{figure*}

The remainder of this paper is organized as follows: Section~\ref{sec: Framework} describes how CSUN's campus electric grid is modeled. In Section~\ref{sec: Model}, we introduce the proposed high-fidelity model for analyzing the impact of EVs on low-voltage distribution systems. Section~\ref{sec: Result} presents the simulation results, and Section~\ref{sec: Conclusion} provides the conclusion.


\section{Derivation of CSUN's electric grid}
\label{sec: Framework}

\mrn{To investigate the impact of EV charging stations on CSUN's electric grid, we developed a detailed model that includes comprehensive building load data, existing EV chargers, cable impedance, and transformer characteristics. This simulation accurately represents the grid's performance under various scenarios, allowing us to assess and address potential challenges from integrating EV charging infrastructure. The derivation of these components is explained as follows:}

\begin{figure}
    \centering
\captionsetup{justification=centering}
\includegraphics[scale=.75,trim= 2cm 18.2cm 11.5cm 1.7cm,clip]{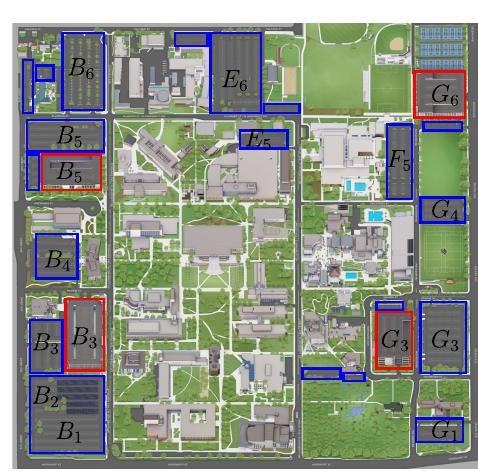}
	 \caption{\mrn{With over 36,000 students, the CSUN campus features various parking areas, including open-air and multi-level structures, to accommodate projected EV charging loads. This high electric load in a compact area stresses the distribution system, requiring load impact analysis to prevent issues like line overloading.}}
	\label{fig:CSUN_MAP}
  \hspace{-0.7cm}
\end{figure}

\begin{figure}
    \centering
    \vspace{-0.8cm}
\captionsetup{justification=centering}
\includegraphics[scale=0.3,trim= 3.5cm 3.0cm 3.cm 0.0cm,clip]{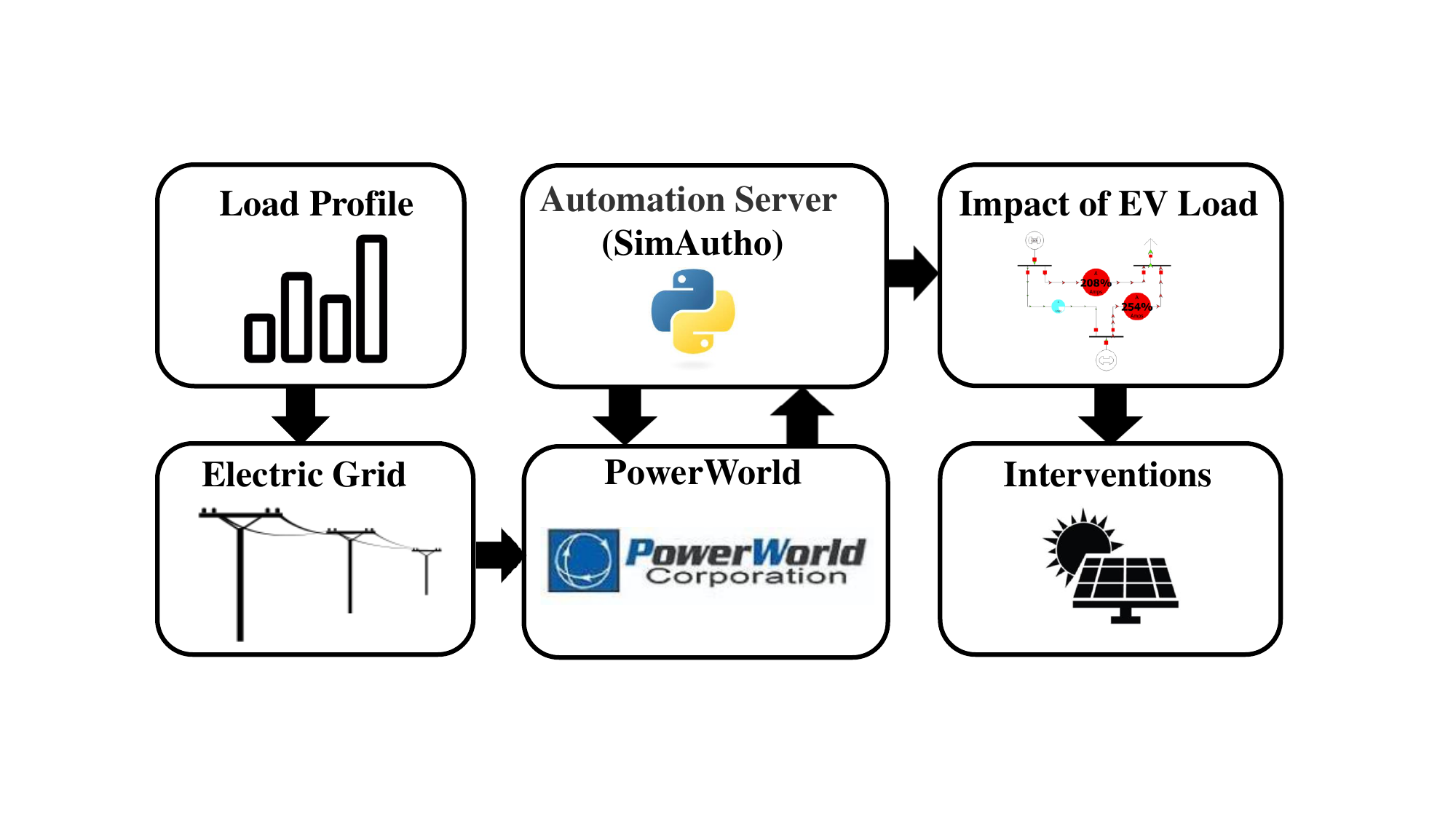}
	 \caption{\mrn{The flowchart of the proposed framework for investigating the impact of EV loads on the power distribution system.}}
	\label{fig:framework}
\end{figure}

\mrn{\textbf{Data Collection and Load Modeling:}}  
\mrn{To simulate the campus grid accurately, we gathered comprehensive load data for all buildings and existing EV chargers by reading meters on each phase of the building cables. Since the system is nearly balanced, a single meter reading per building reliably represented the overall load. This simplified data collection while maintaining accuracy. We also recorded temporal variations in load data to account for daily and seasonal fluctuations, ensuring our model reflects realistic operating conditions.}

\begin{table}
\centering
\begin{scriptsize}
\caption{Parkings' Capacity at CSUN Campus}
\label{Table:Parking_Capacity}
\begin{tabular}{|c|c|c|c|}
\hline 
Parking & Capacity & Parking & Capacity\tabularnewline
\hline 
\hline 
B5  & 420 & B2 & 460\tabularnewline
\hline 
B5 Structure  & 1290 & F2 & 50\tabularnewline
\hline 
E5  & 100 & G1 & 90\tabularnewline
\hline 
F5   & 230 & G3 Structure & 1370\tabularnewline
\hline 
G6 & 50 & G3 & 450\tabularnewline
\hline 
B6 & 460 & B3 Structure & 1760
\tabularnewline
\hline 
E6 & 590 & B4 & 300\tabularnewline
\hline 
G6 Structure & 1300 & G4 & 170\tabularnewline
\hline 
B1 & 480 & B3  & 500\tabularnewline
\hline 
\end{tabular}
\end{scriptsize}
\end{table}
\mrn{\textbf{Cable Resistance and Impedance:}}
\mrn{Determining the resistance and impedance of the cables was crucial for our simulation. We obtained detailed information on the types of cables and the campus's oneline diagram from the Physical Plant Management (PPM) department. The oneline diagram, which includes the coordinates of the cables, allowed us to calculate their lengths accurately. For resistance, calculation, we first identified the specific type of each cable. We then multiplied the ohm/mile value of the corresponding cable type by its length, as derived from the one-line diagram. This provided us with precise resistance values for each segment of the grid. By doing so, we ensured that our model accounts for potential variations in cable performance over time.}

\mrn{\textbf{Transformer Impedance and Capacity:}}
\mrn{Transformers are critical components in the grid, and their accurate modeling is essential for realistic simulations. PPM provided us with detailed specifications for each transformer, including their types, impedance, and capacity. With this data, we were able to incorporate the transformers' characteristics into our grid model, ensuring that their behavior under various load conditions, including the additional demand from EV charging stations, was accurately represented.}

\mrn{\textbf{Distributed Generation (DG):} CSUN has three PV generation sites on campus: ``Parking Lot B2'' with a 467 kW capacity, ``Parking Lot E6'' with a 225 kW capacity, and the ``Student Recreation Center'' with a 1.2 MW capacity. PPM provided us with the PV generation profiles, and we modeled the PV generation as DG in CSUN’s electric grid.}

\section{Proposed Framework}
\label{sec: Model}

\mrn{This paper aims to analyze the impact of various EV loads on the CSUN electrical grid. To achieve this, we have implemented a co-simulation framework, illustrated in Figure \ref{fig:framework}, to evaluate the impacts of EV charging demand on various system parameters such as congestion in lines and transformers. The core analysis of this framework is computed by the PowerWorld simulator. By modeling the proposed electrical network in the PowerWorld simulator, we can comprehensively analyze the system from multiple perspectives. The derivation of these components is explained as follows:}

\begin{figure*}
 \centering
\captionsetup{justification=centering}
\begin{tikzpicture}
\node at (0,0) {\includegraphics[scale=0.9,trim= 4cm 5.26cm 8.08cm 5.245cm  4.5cm,clip]{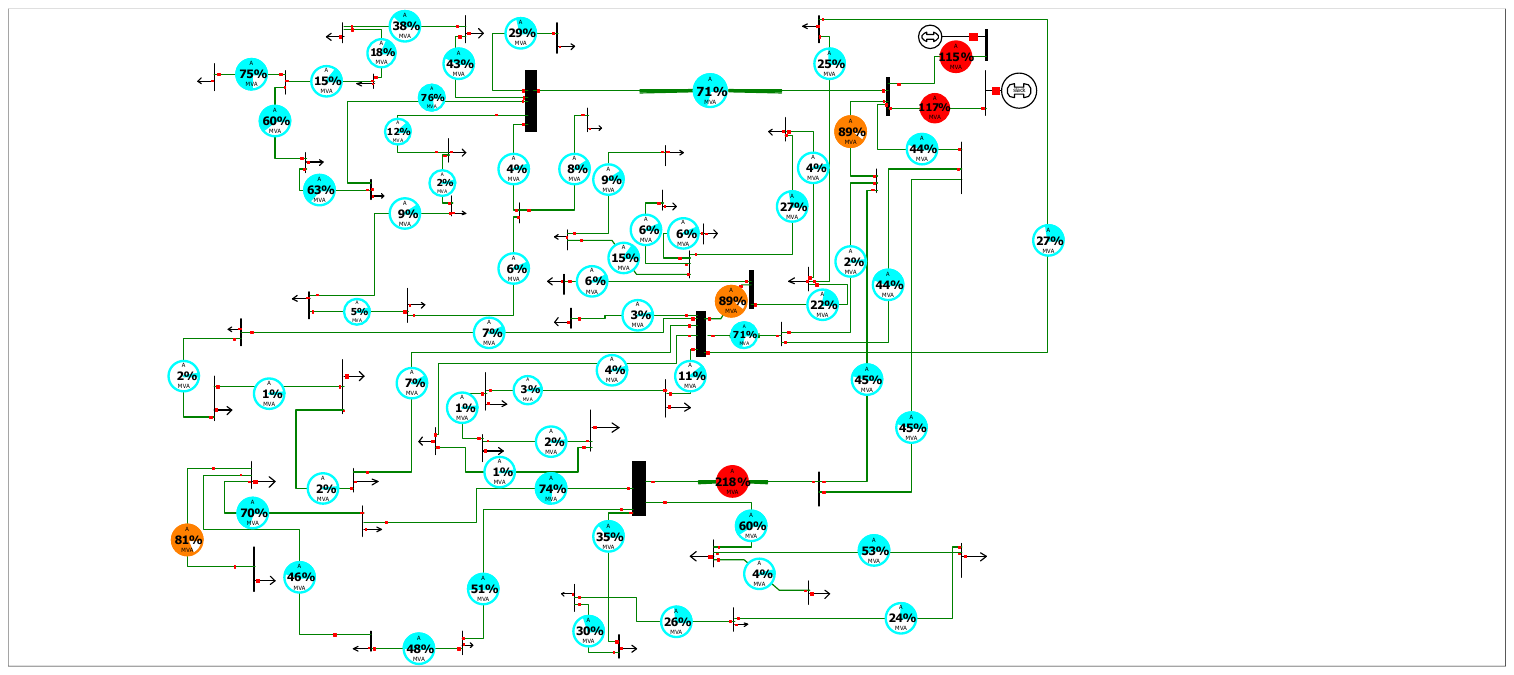}};
 
\node [anchor=west, font=\scriptsize] at (-5.6,5.0) {Physical Plant};
\node [anchor=west, font=\scriptsize] at (-7.6,4.4) {Extended Learning};
\node [anchor=west, font=\scriptsize] at (-4.7,3.7) {Juniper Hall};
\node [anchor=west, font=\scriptsize] at (-3.7,5.1) {Art Design Center};
\node [anchor=west, font=\scriptsize] at (-1.1,4.7) {E6 Mathador Hall};
\node [anchor=west, font=\scriptsize] at (-1.6,3.9) { SubB (LV)};
\node [anchor=west, font=\scriptsize] at (-6.3,2.6) {Education};
\node [anchor=west, font=\scriptsize] at (-4.7,2.5) {Jacaranda Hall};
\node [anchor=west, font=\scriptsize] at (-3.5,3.2) {Redwood Hall};
\node [anchor=west, font=\scriptsize] at (-0.8,3.4) {Sustainability Center};
\node [anchor=west, font=\scriptsize] at (0.4,2.65) {SU(AN)};
\node [anchor=west, font=\scriptsize] at (1.7,3.5) {SU(CP)};
\node [anchor=west, font=\scriptsize] at (2.1,5.1) {Athletic Field};
\node [anchor=west, font=\scriptsize] at (3.1,4.1) {SubB};
\node [anchor=west, font=\scriptsize] at (4.7,3.3) {DWP Pole};
\node [anchor=west, font=\scriptsize] at (4.4,5.0) {Tampa DWP Pole};
\node [anchor=west, font=\scriptsize] at (0.4,2.2) {SU Center};
\node [anchor=west, font=\scriptsize] at (1.0,1.7) {SU (AE)};
\node [anchor=west, font=\scriptsize] at (1.0,1.15) {Student REC};
\node [anchor=west, font=\scriptsize] at (1.8,0.3) {Sub A};
\node [anchor=west, font=\scriptsize] at (0.0,0.5) {Sub A (LV)};
\node [anchor=west, font=\scriptsize] at (-1.7,0.5) {Street LTS};
\node [anchor=west, font=\scriptsize] at (-1.7,1.0) {SU (As)};
\node [anchor=west, font=\scriptsize] at (-1.7,1.8) {SU Admin};
\node [anchor=west, font=\scriptsize] at (-3.7,1.6) {SEQUOIA Hall};
\node [anchor=west, font=\scriptsize] at (-4,0.9) {Oviatt library};
\node [anchor=west, font=\scriptsize] at (-5.6,0.9) {Byramian};
\node [anchor=west, font=\scriptsize] at (-7,0.4) {University Hall};
\node [anchor=west, font=\scriptsize] at (-6.8,-0.4) {Seirra Center};
\node [anchor=west, font=\scriptsize] at (-5.3,-0.2) {Jerome Richfield};
\node [anchor=west, font=\scriptsize] at (-6.9,-1.8) {Parking B2};
\node [anchor=west, font=\scriptsize] at (-4.9,-1.8) {Seirra Hall};
\node [anchor=west, font=\scriptsize] at (-4.9,-3.1) {Manzanita Hall};
\node [anchor=west, font=\scriptsize] at (-6.9,-4) {Parking B3};
\node [anchor=west, font=\scriptsize] at (-5,-5) {Nordhoff Hall};
\node [anchor=west, font=\scriptsize] at (-3.3,-4.3) {Cypress Hall};
\node [anchor=west, font=\scriptsize] at (-1.6,-3.6) {Bookstore};
\node [anchor=west, font=\scriptsize] at (-0.6,-5) {Soraya Hall};
\node [anchor=west, font=\scriptsize] at (0.9,-4.6) {Monterey Hall};
\node [anchor=west, font=\scriptsize] at (-0.1,-3) {Health Center};
\node [anchor=west, font=\scriptsize] at (1.4,-1.8) {Sub C (HV)};
\node [anchor=west, font=\scriptsize] at (-0.4,-1.7) {Sub C};
\node [anchor=west, font=\scriptsize] at (-0.2,-1.3) {Sattelite Plant};
\node [anchor=west, font=\scriptsize] at (-1.6,-1){Chapparal Hall};
\node [anchor=west, font=\scriptsize] at (2.2,-3.6){Chrisholm Hall};
\node [anchor=west, font=\scriptsize] at (5.0,-3.6){Parking G3};
\end{tikzpicture}%
\vspace{-.2cm}
 \caption{CSUN's campus electric grid with regular load and 10\% projected EV loads shows a few line and transformer congestion.}
	\label{fig:powerworld_CSUN_10}
  \vspace{-.3cm}
\end{figure*}

\mrn{First, the elements and load profiles modeled and extracted in Section \ref{sec: Framework} are integrated into the PowerWorld simulator for further analysis. Next, we modeled the load profile of EV charging stations in Python for 15-minute intervals. By normalizing the measured EV load (dividing each interval's load by the maximum load), we obtained coefficients between 0 and 1. These coefficients are then used to update the EV load in PowerWorld by multiplying them with the EV-based load.}
\mrn{We then established a connection between Python and PowerWorld using the Automation Server (SimAutho) to utilize the Application Programming Interface (API) of these two software tools \cite{thayer2020easy}. This connection allows us to update load profile information in the PowerWorld simulator based on the data provided by Python at each step. Based on the updated conditions resulting from the new states of EV load demand, the system performance is analyzed using the PowerWorld simulator. System performance refers to the state of various parameters, such as line congestion.}

\mrn{Depending on the system's condition, PV generation or demand management might be necessary to mitigate EV charging impacts on the grid. These interventions are modeled in Python, and using the API, the updates are communicated to PowerWorld. As shown in Figure \ref{fig:framework}, there is a bidirectional flow of information: load or generation profiles are sent from Python to PowerWorld, and the current system state is read from PowerWorld to Python for decision-making on controlling generation or load demand.}

\mrn{This co-simulation framework allows dynamic analysis and management of the grid's response to varying EV charging demands. By leveraging Python and PowerWorld, we perform detailed simulations and develop strategies to ensure grid stability and efficiency with increasing EV loads. This approach helps understand immediate impacts and provides insights for long-term planning and grid resilience.}


\section{Results and Discussion}
\label{sec: Result}






\begin{figure*}
 \centering
\captionsetup{justification=centering}
\begin{tikzpicture}
\node at (-0.475,-0.0) {\includegraphics[scale=0.93,trim= 5.15cm 5.26cm 8.0cm 5.245cm ,clip]{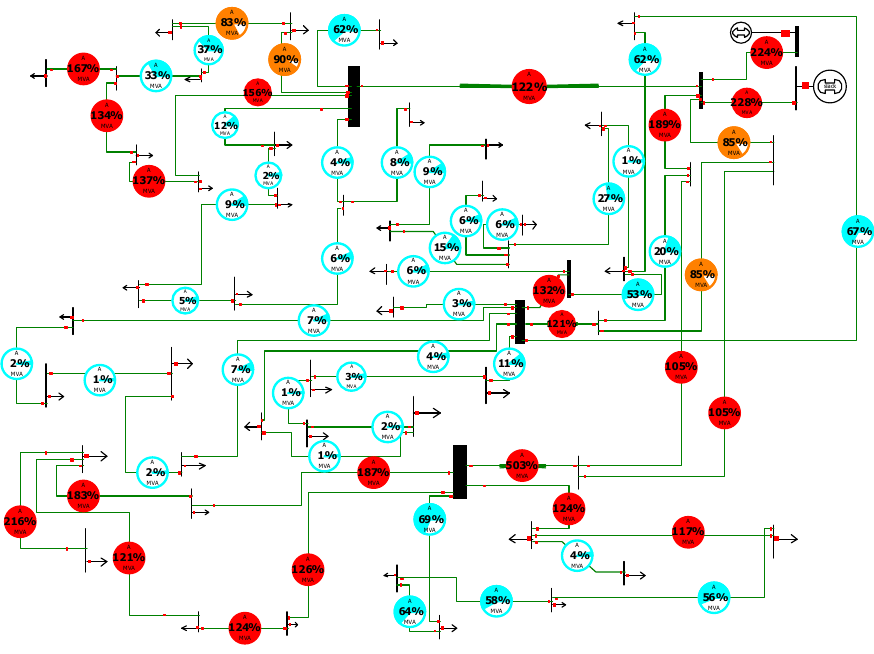}};
 
\node [anchor=west, font=\scriptsize] at (-5.8,5.05) {Physical Plant};
\node [anchor=west, font=\scriptsize] at (-7.8,4.45) {Extended Learning};
\node [anchor=west, font=\scriptsize] at (-4.9,3.75) {Juniper Hall};
\node [anchor=west, font=\scriptsize] at (-3.9,5.15) {Art Design Center};
\node [anchor=west, font=\scriptsize] at (-1.3,4.75) {E6 Mathador Hall};
\node [anchor=west, font=\scriptsize] at (-1.9,3.95) { Sub B (LV)};
\node [anchor=west, font=\scriptsize] at (-6.5,2.65) {Education};
\node [anchor=west, font=\scriptsize] at (-4.9,2.55) {Jacaranda Hall};
\node [anchor=west, font=\scriptsize] at (-3.7,3.25) {Redwood Hall};
\node [anchor=west, font=\scriptsize] at (-1,3.45) {Sustainability Center};
\node [anchor=west, font=\scriptsize] at (0.2,2.7) {SU(AN)};
\node [anchor=west, font=\scriptsize] at (1.5,3.55) {SU(CP)};
\node [anchor=west, font=\scriptsize] at (1.9,5.15) {Athletic Field};
\node [anchor=west, font=\scriptsize] at (2.9,4.15) {Sub B};
\node [anchor=west, font=\scriptsize] at (4.5,3.35) {DWP Pole};
\node [anchor=west, font=\scriptsize] at (4.2,5.05) {Tampa DWP Pole};
\node [anchor=west, font=\scriptsize] at (0.2,2.25) {SU Center};
\node [anchor=west, font=\scriptsize] at (0.8,1.75) {SU (AE)};
\node [anchor=west, font=\scriptsize] at (0.8,1.2) {Student REC};
\node [anchor=west, font=\scriptsize] at (1.6,0.35) {Sub A};
\node [anchor=west, font=\scriptsize] at (-0.2,0.55) {Sub A (LV)};
\node [anchor=west, font=\scriptsize] at (-1.9,0.55) {Street LTS};
\node [anchor=west, font=\scriptsize] at (-1.9,1.05) {SU (As)};
\node [anchor=west, font=\scriptsize] at (-1.9,1.85) {SU Admin};
\node [anchor=west, font=\scriptsize] at (-3.9,1.65) {SEQUOIA Hall};
\node [anchor=west, font=\scriptsize] at (-4.2,0.95) {Oviatt library};
\node [anchor=west, font=\scriptsize] at (-5.8,0.95) {Byramian};
\node [anchor=west, font=\scriptsize] at (-7.2,0.45) {University Hall};
\node [anchor=west, font=\scriptsize] at (-7,-0.35) {Seirra Center};
\node [anchor=west, font=\scriptsize] at (-5.5,-0.15) {Jerome Richfield};
\node [anchor=west, font=\scriptsize] at (-7.1,-1.75) {Parking B2};
\node [anchor=west, font=\scriptsize] at (-5.1,-1.75) {Seirra Hall};
\node [anchor=west, font=\scriptsize] at (-5.1,-3.05) {Manzanita Hall};
\node [anchor=west, font=\scriptsize] at (-7.1,-3.95) {Parking B3};
\node [anchor=west, font=\scriptsize] at (-5.2,-4.95) {Nordhoff Hall};
\node [anchor=west, font=\scriptsize] at (-3.5,-4.25) {Cypress Hall};
\node [anchor=west, font=\scriptsize] at (-1.8,-3.55) {Bookstore};
\node [anchor=west, font=\scriptsize] at (-0.8,-4.95) {Soraya Hall};
\node [anchor=west, font=\scriptsize] at (0.7,-4.55) {Monterey Hall};
\node [anchor=west, font=\scriptsize] at (-0.3,-2.95) {Health Center};
\node [anchor=west, font=\scriptsize] at (1.2,-1.75) {Sub c (HV)};
\node [anchor=west, font=\scriptsize] at (-0.6,-1.65) {Sub c};
\node [anchor=west, font=\scriptsize] at (-0.4,-1.25) {Sattelite Plant};
\node [anchor=west, font=\scriptsize] at (-1.8,-0.95){Chapparal Hall};
\node [anchor=west, font=\scriptsize] at (2.0,-3.55){Chrisholm Hall};
\node [anchor=west, font=\scriptsize] at (4.8,-3.55){Parking G3};
\end{tikzpicture}%
\vspace{-.3cm}
 \caption{CSUN's campus electric grid, with regular load and 25\% projected EV loads, shows significant line and transformer congestion, primarily along the path connecting the LADWP service point to the campus and near parking lots with high EV loads.}
	\label{fig:powerworld_CSUN_25}
 \vspace{-.3cm}
\end{figure*}

\mrn{In this section, we apply the proposed framework to the system described in Section \ref{sec: Framework}. Using a predefined load profile, we analyze the impacts of EV charging on CSUN's electric network through three scenarios at 9:00 AM, when most EVs are connected. The case studies consider the impact of EV loads and interventions like load management and renewable energy resources to alleviate the grid's EV load impact.}


\mrn{\textbf{Scenario I. Impact of EV Charging:}
In this scenario, we examine the impacts of EV charging without PV generation or demand management to understand how the grid handles increased load solely from EV charging stations. By analyzing parameters like line congestion, we identify potential stress points in the system. Fig.~\ref{fig:powerworld_CSUN} shows the simulated CSUN power grid without EV loads, revealing no line congestion. Figs.~\ref{fig:powerworld_CSUN_10} and~\ref{fig:powerworld_CSUN_25} illustrate the impact of 10\% and 25\% EV penetration, respectively. At 10\% EV penetration, two lines and one transformer are congested, while 25\% penetration causes extreme congestion in 18 lines and three transformers. The X\% EV load penetration is calculated by multiplying X\% of the parking capacity by the power drawn by each charging station, which ranges from 7-19 kW.
For this paper, we use 10 kW per charging station. These congested elements are between large parking lots and the entry points from the Los Angeles Department of Water and Power's service points. Preventive actions like load management or distributed generation are necessary to alleviate these issues, as discussed in the following sections.}

\mrn{\textbf{Scenario II. Impact of Distributed Generation:}
We assess the impact of integrating photovoltaic (PV) systems into the grid model to offset the additional load from EV chargers. This scenario evaluates how renewable energy can mitigate the negative effects of increased EV load, especially during peak times. 
Incorporating photovoltaic (PV) generation into the CSUN grid with 10\% electric vehicle (EV) loads reduces the number of congested lines from six to two, as shown in Fig.\ref{fig:powerworld_CSUN_10}. Additionally, the severity of congestion in these two lines is lower than in their counterparts in Fig.\ref{fig:powerworld_CSUN_10}. Therefore, adding PV generation can alleviate line congestion caused by EV loads. However, some lines remain congested or near capacity, indicating that additional interventions like load shifting may be necessary.}


\begin{table}
\begin{scriptsize}
    \centering
\caption{Distribution of Line Flow Percentages Across Different Scenarios}
\begin{tabular}{|s|c|c|c|c|}
\hline 
 & \multicolumn{4}{c|}{\thead{Power Flow Percentage intervals\\ for Lines and Transformers}}\tabularnewline
\hline 
\hline 
\thead{Scenarios} & \thead{40--80} & \thead{80--100} & \thead{100--150} & \thead{\textgreater 150}\tabularnewline
\hline 
\thead{Current Grid} & \thead{2} & {\thead{--}} & {\thead{--}} & {\thead{--}}\tabularnewline
\hline 
\thead{Current Grid \& 10\%\\ EV Load} & \thead{16} & \thead{3} & \thead{2} & {\thead{--}}\tabularnewline
\hline 
\thead{Current Grid with 25\%\\ EV Load} & \thead{9} & \thead{4} & \thead{10} & \thead{8}\tabularnewline
\hline 
\thead{Current Grid \& 25\%\\ EV Load \& PV Generation} & \thead{16} & \thead{2} & \thead{1} & \thead{1}\tabularnewline
\hline 
\thead{Current Grid \& 25\%\\ EV Load \& PV Generation\\ \& Load Management} & \thead{11} & \thead{3} & \thead{0} & \thead{1}\tabularnewline
\hline 
\end{tabular}
    \label{tab:my_label}
    \end{scriptsize}
    \vspace{-0.6cm}
\end{table}

\begin{figure*}
 \centering
\captionsetup{justification=centering}
\begin{tikzpicture}
\node at (-0.58,0.1) {\includegraphics[scale=0.92,trim= 4.55cm 5.30cm 8.1cm 5.25cm  4.5cm,clip]{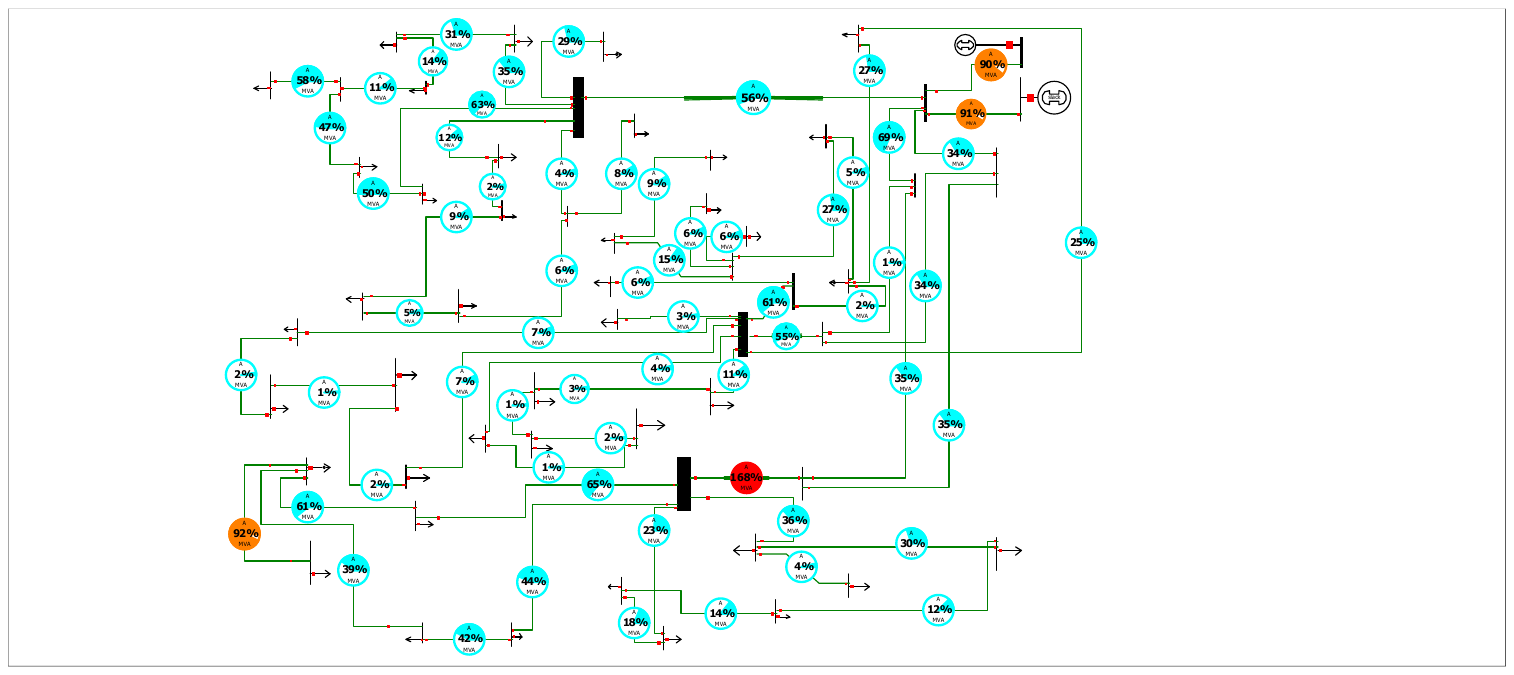}};
 
\node [anchor=west, font=\scriptsize] at (-5.8,5.05) {Physical Plant};
\node [anchor=west, font=\scriptsize] at (-7.8,4.45) {Extended Learning};
\node [anchor=west, font=\scriptsize] at (-4.9,3.75) {Juniper Hall};
\node [anchor=west, font=\scriptsize] at (-3.9,5.15) {Art Design Center};
\node [anchor=west, font=\scriptsize] at (-1.3,4.75) {E6 Mathador Hall};
\node [anchor=west, font=\scriptsize] at (-1.9,3.95) { Sub B (LV)};
\node [anchor=west, font=\scriptsize] at (-6.5,2.65) {Education};
\node [anchor=west, font=\scriptsize] at (-4.9,2.55) {Jacaranda Hall};
\node [anchor=west, font=\scriptsize] at (-3.7,3.25) {Redwood Hall};
\node [anchor=west, font=\scriptsize] at (-1,3.45) {Sustainability Center};
\node [anchor=west, font=\scriptsize] at (0.2,2.7) {SU(AN)};
\node [anchor=west, font=\scriptsize] at (1.5,3.55) {SU(CP)};
\node [anchor=west, font=\scriptsize] at (1.9,5.15) {Athletic Field};
\node [anchor=west, font=\scriptsize] at (2.9,4.15) {Sub B};
\node [anchor=west, font=\scriptsize] at (4.5,3.35) {DWP Pole};
\node [anchor=west, font=\scriptsize] at (4.2,5.05) {Tampa DWP Pole};
\node [anchor=west, font=\scriptsize] at (0.2,2.25) {SU Center};
\node [anchor=west, font=\scriptsize] at (0.8,1.75) {SU (AE)};
\node [anchor=west, font=\scriptsize] at (0.8,1.2) {Student REC};
\node [anchor=west, font=\scriptsize] at (1.6,0.35) {Sub A};
\node [anchor=west, font=\scriptsize] at (-0.2,0.55) {Sub A (LV)};
\node [anchor=west, font=\scriptsize] at (-1.9,0.55) {Street LTS};
\node [anchor=west, font=\scriptsize] at (-1.9,1.05) {SU (As)};
\node [anchor=west, font=\scriptsize] at (-1.9,1.85) {SU Admin};
\node [anchor=west, font=\scriptsize] at (-3.9,1.65) {SEQUOIA Hall};
\node [anchor=west, font=\scriptsize] at (-4.2,0.95) {Oviatt library};
\node [anchor=west, font=\scriptsize] at (-5.8,0.95) {Byramian};
\node [anchor=west, font=\scriptsize] at (-7.2,0.45) {University Hall};
\node [anchor=west, font=\scriptsize] at (-7,-0.35) {Seirra Center};
\node [anchor=west, font=\scriptsize] at (-5.5,-0.15) {Jerome Richfield};
\node [anchor=west, font=\scriptsize] at (-7.1,-1.75) {Parking B2};
\node [anchor=west, font=\scriptsize] at (-5.1,-1.75) {Seirra Hall};
\node [anchor=west, font=\scriptsize] at (-5.1,-3.05) {Manzanita Hall};
\node [anchor=west, font=\scriptsize] at (-7.1,-3.95) {Parking B3};
\node [anchor=west, font=\scriptsize] at (-5.2,-4.95) {Nordhoff Hall};
\node [anchor=west, font=\scriptsize] at (-3.5,-4.25) {Cypress Hall};
\node [anchor=west, font=\scriptsize] at (-1.8,-3.55) {Bookstore};
\node [anchor=west, font=\scriptsize] at (-0.8,-4.95) {Soraya Hall};
\node [anchor=west, font=\scriptsize] at (0.7,-4.55) {Monterey Hall};
\node [anchor=west, font=\scriptsize] at (-0.3,-2.95) {Health Center};
\node [anchor=west, font=\scriptsize] at (1.2,-1.75) {Sub c (HV)};
\node [anchor=west, font=\scriptsize] at (-0.6,-1.65) {Sub c};
\node [anchor=west, font=\scriptsize] at (-0.4,-1.25) {Sattelite Plant};
\node [anchor=west, font=\scriptsize] at (-1.8,-0.95){Chapparal Hall};
\node [anchor=west, font=\scriptsize] at (2.0,-3.55){Chrisholm Hall};
\node [anchor=west, font=\scriptsize] at (4.8,-3.55){Parking G3};
\end{tikzpicture}%
\vspace{-.3cm}
 \caption{CSUN's campus electric grid with regular load, 25\% projected EV loads, and three photovoltaic distributed generation units on Buses ``Parking B2'', ``Student REC'', and ``E6 Matador Hall''. In this case, an EV load management system is considered to prevent all EV loads from being connected simultaneously and to shift a portion of these loads to future time intervals. Incorporating the distributed generation and applying EV load management drastically improves line congestion compared to Fig.~\ref{fig:powerworld_CSUN_25}.}
	\label{fig:powerworld_CSUN_25_PV_LM}
 \vspace{-.3cm}
\end{figure*}

\mrn{\textbf{Scenario III. Impacts of PV Generation and Demand Management}
In this case, demand management strategies, such as shifting EV charging times and implementing controlled schedules, are applied alongside PV generation to mitigate grid impact. Fig.~\ref{fig:powerworld_CSUN_25_PV_LM} shows CSUN's grid with 25\% projected EV loads, PV generation at three sites, and an EV load management mechanism that connects only one-third of EV loads at a time. This combined approach significantly reduces line congestion compared to Fig.~\ref{fig:powerworld_CSUN_25}, indicating that PV generation and load management together can effectively support the integration of large EV loads into the grid.}

\mrn{Table~\ref{tab:my_label} shows the distribution of line flow percentages across various scenarios. Adding EV loads increases congestion, as indicated by more lines with power flows exceeding their capacity. For instance, the current grid has 2 lines with 40-80\% flow, but a 10\% EV load raises this to 16 lines. A 25\% EV load results in 10 lines with 100-150\% flow and 8 lines exceeding 150\%. However, integrating photovoltaic (PV) generation alleviates this congestion, reducing overloaded lines. Implementing load management further optimizes power flow distribution, enhancing grid efficiency. This combined approach is crucial for mitigating the impact of EV loads.}

\section{Conclusion}
\label{sec: Conclusion}
\mrn{This paper presents a detailed and realistic model of CSUN's electric grid for various electrical analyses. We explored scenarios to analyze the stress of EV charging stations on the grid and how it can be alleviated using PV generation and load management. Each scenario helped us understand the distinct and combined effects of EV charging, PV generation, and demand management on CSUN's grid performance.
By comparing these scenarios, we identified effective strategies for integrating EV charging infrastructure while maintaining grid reliability and efficiency. Our analyses highlight the benefits and challenges of each approach, offering insights into potential solutions for managing the growing demand for EV charging at CSUN.
Additionally, we performed load flow analyses under different scenarios, such as peak and off-peak hours, and assessed the impact of various EV charger deployment configurations on transformer loading. This allowed us to identify potential bottlenecks and stress points, guiding necessary infrastructure upgrades.
Overall, this study provides valuable information for planning and implementing EV charging infrastructure, ensuring grid stability and efficiency while accommodating the increasing demand for electric vehicles. This research is ongoing to determine the optimal locations for charging stations on CSUN's campus.}
\vspace{-.2cm}


\section*{Acknowledgement}
This research is developed as a part of the project ``Climate Action - Community-driven eLectric vEhicle
chArging solutioN (CA-CLEAN),'' funded by the Climate Action University of California. It is also supported by the US Department of Education grant, SECURE for Student Success ($\textit{SfS}^2$).
\vspace{-.2cm}


\bibliographystyle{IEEEtran}
\IEEEtriggeratref{20}
\bibliography{ref}

\begin{thebibliography}{10}
\providecommand{\url}[1]{#1}
\csname url@samestyle\endcsname
\providecommand{\newblock}{\relax}
\providecommand{\bibinfo}[2]{#2}
\providecommand{\BIBentrySTDinterwordspacing}{\spaceskip=0pt\relax}
\providecommand{\BIBentryALTinterwordstretchfactor}{4}
\providecommand{\BIBentryALTinterwordspacing}{\spaceskip=\fontdimen2\font plus
\BIBentryALTinterwordstretchfactor\fontdimen3\font minus \fontdimen4\font\relax}
\providecommand{\BIBforeignlanguage}[2]{{%
\expandafter\ifx\csname l@#1\endcsname\relax
\typeout{** WARNING: IEEEtran.bst: No hyphenation pattern has been}%
\typeout{** loaded for the language `#1'. Using the pattern for}%
\typeout{** the default language instead.}%
\else
\language=\csname l@#1\endcsname
\fi
#2}}
\providecommand{\BIBdecl}{\relax}
\BIBdecl

\bibitem{muratori2021rise}
M.~Muratori, M.~Alexander, D.~Arent, M.~Bazilian, P.~Cazzola, E.~M. Dede, J.~Farrell, C.~Gearhart, D.~Greene, A.~Jenn \emph{et~al.}, ``The rise of electric vehicles—2020 status and future expectations,'' \emph{Progress in Energy}, vol.~3, no.~2, p. 022002, 2021.

\bibitem{secchi2023smart}
M.~Secchi, G.~Barchi, D.~Macii, and D.~Petri, ``Smart electric vehicles charging with centralised vehicle-to-grid capability for net-load variance minimisation under increasing ev and pv penetration levels,'' \emph{Sustainable Energy, Grids and Networks}, vol.~35, p. 101120, 2023.

\bibitem{crozier2020opportunity}
C.~Crozier, T.~Morstyn, and M.~McCulloch, ``The opportunity for smart charging to mitigate the impact of electric vehicles on transmission and distribution systems,'' \emph{Applied Energy}, vol. 268, p. 114973, 2020.

\bibitem{yong2023electric}
J.~Y. Yong, W.~S. Tan, M.~Khorasany, and R.~Razzaghi, ``Electric vehicles destination charging: An overview of charging tariffs, business models and coordination strategies,'' \emph{Renewable and Sustainable Energy Reviews}, vol. 184, p. 113534, 2023.

\bibitem{narimani2019efficient}
M.~R. Narimani, A.~Azizivahed, and E.~Naderi, ``An efficient scenario-based stochastic energy management of distribution networks with distributed generation, pv module, and energy storage,'' \emph{arXiv preprint arXiv:1910.07109}, 2019.

\bibitem{ramadan2018smart}
H.~Ramadan, A.~Ali, M.~Nour, and C.~Farkas, ``Smart charging and discharging of plug-in electric vehicles for peak shaving and valley filling of the grid power,'' in \emph{2018 Twentieth International Middle East Power Systems Conference (MEPCON)}.\hskip 1em plus 0.5em minus 0.4em\relax IEEE, 2018, pp. 735--739.

\bibitem{nafisi2015two}
H.~Nafisi, S.~M.~M. Agah, H.~A. Abyaneh, and M.~Abedi, ``Two-stage optimization method for energy loss minimization in microgrid based on smart power management scheme of phevs,'' \emph{IEEE Transactions on Smart Grid}, vol.~7, no.~3, pp. 1268--1276, 2015.

\bibitem{clairand2017tariff}
J.-M. Clairand, J.~R. Garc{\'\i}a, C.~{\'A}. Bel, and P.~P. Sarmiento, ``A tariff system for electric vehicle smart charging to increase renewable energy sources use,'' in \emph{2017 IEEE PES Innovative Smart Grid Technologies Conference-Latin America (ISGT Latin America)}.\hskip 1em plus 0.5em minus 0.4em\relax IEEE, 2017, pp. 1--6.

\bibitem{asghari2019method}
B.~Asghari, M.~R. Narimani, and R.~Sharma, ``Method for operation of energy storage systems to reduce demand charges and increase photovoltaic (pv) utilization,'' Jan.~31 2019, uS Patent App. 16/006,239.

\bibitem{yuvaraj2024comprehensive}
T.~Yuvaraj, K.~Devabalaji, J.~A. Kumar, S.~B. Thanikanti, and N.~Nwulu, ``A comprehensive review and analysis of the allocation of electric vehicle charging stations in distribution networks,'' \emph{IEEE Access}, 2024.

\bibitem{nour2019smart}
M.~Nour, S.~M. Said, A.~Ali, and C.~Farkas, ``Smart charging of electric vehicles according to electricity price,'' in \emph{2019 international conference on innovative trends in computer engineering (ITCE)}.\hskip 1em plus 0.5em minus 0.4em\relax IEEE, 2019, pp. 432--437.

\bibitem{narimani2019demand}
M.~R. Narimani, ``Demand side management for homes in smart grids,'' in \emph{2019 North American Power Symposium (NAPS)}.\hskip 1em plus 0.5em minus 0.4em\relax IEEE, 2019, pp. 1--6.

\bibitem{narimani2017multi}
M.~R. Narimani, Maigha, J.-Y. Joo, and M.~Crow, ``Multi-objective dynamic economic dispatch with demand side management of residential loads and electric vehicles,'' \emph{Energies}, vol.~10, no.~5, p. 624, 2017.

\bibitem{narimani2015dynamic}
M.~R. Narimani, J.-Y. Joo, and M.~L. Crow, ``Dynamic economic dispatch with demand side management of individual residential loads,'' in \emph{2015 North American Power Symposium (NAPS)}.\hskip 1em plus 0.5em minus 0.4em\relax IEEE, 2015, pp. 1--6.

\bibitem{azizivahed2017new}
A.~Azizivahed, E.~Naderi, H.~Narimani, M.~Fathi, and M.~R. Narimani, ``A new bi-objective approach to energy management in distribution networks with energy storage systems,'' \emph{IEEE Transactions on Sustainable Energy}, vol.~9, no.~1, pp. 56--64, 2017.

\bibitem{rahman2022comprehensive}
S.~Rahman, I.~A. Khan, A.~A. Khan, A.~Mallik, and M.~F. Nadeem, ``Comprehensive review \& impact analysis of integrating projected electric vehicle charging load to the existing low voltage distribution system,'' \emph{Renewable and Sustainable Energy Reviews}, vol. 153, p. 111756, 2022.

\bibitem{boyaci2022spatio}
O.~Boyaci, M.~R. Narimani, K.~Davis, and E.~Serpedin, ``Spatio-temporal failure propagation in cyber-physical power systems,'' in \emph{2022 3rd International Conference on Smart Grid and Renewable Energy (SGRE)}.\hskip 1em plus 0.5em minus 0.4em\relax IEEE, 2022, pp. 1--6.

\bibitem{narimani2016reliability}
M.~R. Narimani, P.~J. Nauert, J.-Y. Joo, and M.~L. Crow, ``Reliability assesment of power system at the presence of demand side management,'' in \emph{2016 IEEE Power and Energy Conference at Illinois (PECI)}.\hskip 1em plus 0.5em minus 0.4em\relax IEEE, 2016, pp. 1--5.

\bibitem{thayer2020easy}
B.~L. Thayer, Z.~Mao, Y.~Liu, K.~Davis, and T.~Overbye, ``Easy simauto (esa): A python package that simplifies interacting with powerworld simulator,'' \emph{Journal of Open Source Software}, vol.~5, no.~50, p. 2289, 2020.

\end{thebibliography}
\end{document}